\documentclass[leqno]{amsart}
\usepackage{amsthm}
\usepackage{amsfonts}
\usepackage{amsmath}
\usepackage{latexsym,euscript}
\usepackage{graphicx}
\begin{document}

\newcommand{\RR}{\mathbb{R}}
\newcommand{\CC}{\mathbb{C}} 
\newcommand{\ZZ}{\mathbb{Z}}
\newcommand{\NN}{\mathbb{N}}
\newcommand{\FF}{\mathbb{F}}
\newcommand{\I}{\mathrm{i}}
\newcommand{\E}{\mathrm{e}}
\newcommand{\D}{\,\mathrm{d}}
\def\bin#1#2{\biggl({#1\atop#2}\biggr)}
\def\dbin#1#2{{\displaystyle\biggl({#1\atop#2}\biggr)}}
\def\Dbin#1#2{{\displaystyle\Biggl({#1\atop#2}\Biggr)}}
\def\Cale#1{\EuScript #1}
\newcommand{\PP}{{\Cale P}}
\newcommand{\LL}{{\Cale L}}
\newcommand{\DS}{\displaystyle}

\oddsidemargin 16.5mm
\evensidemargin 16.5mm

\thispagestyle{plain}




\vspace{5cc}
\begin{center}

{\large\bf  OSTROWSKI TYPE INEQUALITIES AND SOME SELECTED QUADRATURE FORMULAE
\rule{0mm}{6mm}\renewcommand{\thefootnote}{}
\footnotetext{ 

\rule{2.4mm}{0mm}Keywords and Phrases. Inequalities of Ostrowski type, norm, quadrature rules, Peano kernel, best constant.
}}

\vspace{1cc}
{\large\it Gradimir V. Milovanovi\'c}

\vspace{1cc}
{{\it 
Dedicated to the Memory of  Professor Dragoslav S. Mitrinovi\'c (1908--1995)}}

\vspace{1cc}
\parbox{24cc}{{\small
Abstract. Some selected Ostrowski type inequalities and  a connection with numerical integration are studied in this survey paper, which is dedicated to the memory of Professor D.\ S.\ Mitrinovi\'c, who left us 25 years ago.   His significant influence to the development of the theory of inequalities is briefly given in the first section of this paper. Beside some basic facts on quadrature formulas  and an approach for estimating the error term using Ostrowski type inequalities and Peano kernel techniques, we  give several examples of selected  quadrature formulas and the corresponding inequalities,  including the basic Ostrowski's inequality (1938),   inequality of Milovanovi\'c and Pe\v cari\'c (1976) and its modifications, inequality of Dragomir, Cerone and Roumeliotis (2000),  symmetric inequality of Guessab and Schmeisser (2002) and asymmetric inequality of Franji\'c (2009), as well as four point symmetric inequalites by Alomari (2012) and  a variant with double internal nodes given by Liu and Park (2017).}}
\end{center}

\vspace{1.5cc}
\begin{center}
{\bf 1. MITRINOVI\'C'S INFLUENCE TO THE THEORY OF INEQUALITIES}
\end{center}

My university and scientific career began in the seventies of the last century and is related to Professor Dragoslav S. Mitrinovi\'c
(1908--1995), who at that time was the Head of the Department of Mathematics at the Faculty of Electrical Engineering in Belgrade and founder and  Editor-in-Chief of the journal {\sc Univ.\  Beograd.\ Publ.\ Elektrotehn.\ Fak.\ Ser.\ Mat.\ Fiz.}, started in 1956,  that continue to live through today's {\sc Appl.\ Anal.\  Discrete Math.}  journal (name was changed  in 2007). Although Mitrinovi\'c  dealt with differential and functional equations, as a prominent member of the Belgrade School of Mathematics founded by Mihailo Petrovi\'c Alas (1868--1943) \cite{mpa}, as well as  other areas of real and complex analysis, special functions, number theory, etc., but the inequalities were his greatest passion in mathematics. He was involved in all kinds of inequalities. He often used to say  {\it``There are no equalities, even in the human life, the inequalities are always met.''} (for details see \cite{GVM-1998}, as well as the complete book \cite{RPI}). What should be emphasized is that Mitrinovi\'c was a scientist who always advocated scientific honesty. 
He always warned his associates that they must correctly cite the results of other authors, no matter what personal relationship they have with them. 

 Although I published the first few scientific papers in the field of numerical analysis (especially in iterative processes) and functional equations, Mitrinovi\'c's influence prevailed and I began to deal with the theory of inequalities. His famous monograph ``{\it Analytic Inequalities}'' \cite{AI70} published in 1970 by Springer was at that time an extraordinary inspiration not only for me, but also for many in the world, especially mathematicians of the younger generation, who were able to find interesting topics and sources for their research there. 
 In their review of this monograph R.\ A.\ Askey and R.\ P.\ Boas Jr.\ in {\sc Math.\ Reviews} (MR 274686 (43 \#448)), compared it with previous famous monographs on inequalities written by  G.\ H.\ Hardy, J.\ E.\ Littlewood and G.\ P\'olya \cite{HLP}  and  by E.\ F.\ Beckenbach and R.\ Bellman \cite{BB}, said that 
``{\it Anyone interested in the subject will have to have all three: Hardy-Littlewood-P\'olya for its exhaustive treatment of the classical inequalities and for its thorough discussion of advanced topics that do not appear in other books $\ldots$;  Beckenbach-Bellman for its wide range both of methods and of topics; and Mitrinovi\'c for topics that are in neither of the other books; for its thorough bibliographies; and for an extensive collection of special inequalities, many of which are not otherwise easily accessible, and some of which appear here for the first time. By searching the literature the author has recovered many interesting inequalities that would otherwise have been forgotten. Although what appeals to one analyst need not appeal to another, almost anybody is sure to find interesting things in this book.}''  Describing this three-part book by Mitrinovi\'c, they emphasize that ``{\it The third and most significant part of the book contains some 450 (by the author's count) particular inequalities, loosely arranged according to subject matter. This is a valuable source of material and many of the inequalities could serve as starting points for more general theories.}''

This is exactly what happened in the following period. These inequalities have attracted the attention of many authors, leading to rapid progress and the creation of a theory of inequality based on the linking of many particular inequalities and their generalizations. My first interest was also related to the third part of this monograph, precisely to Section  3.7  on the so-called  {\it Integral Inequalities} \cite[pp.~289--309]{AI70}, as well as to {\it Miscellaneous Inequalities}, which are given later in Section 3.9.  I was particularly drawn to those integral inequalities, given  with appropriate references, as items {\bf3.7.22} (Mackey \cite{Ma1947}), {\bf3.7.23} (Ostrowski \cite{Os1938}), {\bf3.7.24} (Iyengar \cite{Iy1938}), and {\bf3.7.29} (Zmorovi\v c \cite{Zm1956}), and working on them I obtained several extensions and generalizations \cite{DjGVM1976,GM2,GM5}, which I 
 used in order to estimate the error terms in some general quadrature formulas, among other things (see also a later result on Iyengar inequality obtained jointly with Pe\v cari\'c \cite{GM4Y}). 
 
 For example, using Ostrowski's inequality   
\cite{Os1938}
 \begin{equation}\label{Os1}
\Biggl\vert \ f(x)-\frac{1}{b-a}\int_{a}^{b}f(t){\D}t\Biggr\vert \leq  
\Biggl[ \frac{1}{4}+\frac{\left(x-\frac{a+b}{2}\right)^{2}}{(
b-a)^2}\Biggr] (b-a)\left\Vert f'\right\Vert _{\infty },\quad x\in[a,b],
\end{equation}
which holds for each continuous function $f:[a,b] \to\RR$,  differentiable on $(a,b)$, with bounded derivative 
$$\|f'\|_{\infty}=\sup\limits_{x\in[a,b]}|f'(t)|<\infty,$$ 
for such functions on $[0,1]$ the following simple estimate   \cite{GM2} 
\[ \left\vert \int_0^1f(x){\D}x-\frac{1}{n}\sum_{k=1}^n\lambda_kf(x_k)\right\vert \le \frac{M}{2}\sum_{k=1}^n \left[(x_k-a_{k-1})^2+(a_k-x_k)^2\right]\]
was proved, where $0=a_0<a_1<a_2<\cdots<a_n=1$ and
\[\lambda_k=a_k-a_{k-1},\quad a_{k-1}\le x_k\le a_k\ \ (k=1,\ldots,n).\]
In the same paper \cite{GM2} we also  proved the multidimensional version 
of (\ref{Os1}), including the weighted case, as well as the corresponding applications in numerical  integration over the domain $D=\bigl\{(x_1,\ldots ,x_m\}\,:\,a_i<x_i<b_i\ (i=1,\ldots,m)\bigr\}$.\\

\noindent
{\bf Theorem 1.1.} {\it
Let $f:\RR^m\to\RR$ be  a differentiable
function defined on $\overline{D}$ and let $\left\vert \dfrac{\partial f}{ 
\partial x_i}\right\vert \leq M_i$ $(M_i>0;\ i=1,\ldots ,m)$ in $D$. 
Then, for every $\left( x_{1},\ldots
,x_{m}\right) \in $ $\overline{D}$,
\begin{align}\label{A-2} 
\Biggl\vert f(x_1,\ldots ,x_m) -\frac{1}{\prod\limits_
{i=1}^m (b_i-a_i)}&\int_{a_1}^{b_1} 
\cdots\int_{a_m}^{b_m}f(y_{1},\ldots
,y_{m}) {\D}y_{1}\cdots {\D}y_{m}\Biggr\vert \\ 
&\le \sum_{i=1}^m\left[ \dfrac{1}{4}+\frac{\left(
x_{i}-\dfrac{a_{i}+b_{i}}{2}\right) ^{2}}{\left( b_{i}-a_{i}\right) ^{2}}%
\right] \left( b_{i}-a_{i}\right) M_{i}.  
\notag
\end{align}}

 It seems that the previous results were the first application of Ostrowski's inequality in numerical integration for getting  estimates of the remainder term in  composite quadrature formulas.  
 Also, we used  {\bf3.9.71}, i.e., the Landau inequality  
 $|f'(x)|\le 2$ $(x\in I)$, which holds for all real functions $x\mapsto f(x)$ on an interval $I$,  of length not less than $2$, for which $|f(x)|\le 1$ and $|f''(x)|\le 1$ (see \cite{La1914}), as well as its  generalization 
 \begin{equation}\label{AvAlj}
 \bigl|\varphi'(x)-\varphi(1)+\varphi(0\bigr|
         \le\frac{1}{2}-x+x^2,\quad 0\le x\le 1,	
 \end{equation}
proved  by Avakumovi\'c and Aljan\v ci\'c \cite{AvAlj1950}, by geometric arguments
under the condition  $|\varphi''(x)|\le 1$ for $0\le x\le 1$. Here, the polynomial $x\mapsto \frac{1}{2}-x+x^2$ is the best possible, as well as the constant $2$ in the Landau inequality. 

Otherwise, there are several generalizations of the Landau  
result in many senses. Our generalization was related to twice Fr\'chet-differentiable operators  $F:X\to Y$, where $X$ and $Y$  are Banach spaces (see \cite{GM3,GM1}).  

The inequality (\ref{AvAlj}) is connected with the  Ostrowski inequality (\ref{Os1}) and it can be seen if we take  
\[\varphi(x)=\frac{1}{M(b-a)}\int_0^x f(a+(b-a)t){\D}t.\]
Then   (\ref{AvAlj}) reduces to  (\ref{Os1}). 

Several monographs have been also appeared after Mitrinovi\'c's monograph \cite{AI70}. We mention here only a few of them: {\it Means and Their Inequalities}  \cite{BMV88} by Bullen, Mitrinovi\'c and Vas\'c, {\it Inequalities Involving Functions and Their Integrals and Derivatives} \cite{MPF91} and  {\it Classical and New Inequalities in Analysis} \cite{MPF93} by Mitrinovi\'c, Pe\v cari\'c and Fink, and  {\it Topics in Polynomials: Extremal Problems, Inequalities, Zeros} \cite{MMR94} by Milovanovi\'c, Mitrinovi\'c and Rassias.

From today's point of view, we can notice that after the mentioned period and \cite[Chp.\ XV]{MPF91}, Ostrowski's inequality (\ref{Os1}) became a challenge for many researchers, so according to {\it Math.\ Review} and to Dragomir's survey paper \cite{Drag17}, there are a few hundreds of published papers with the phrase ``{\it Ostrowski type of inequality}'' in the title, and even an edited book by Dragomir and Rassias \cite{DragRas}, as well as a nice monograph by Franji\'c, Pe\v cari\'c, Peri\'c and Vukeli\'c \cite{Mon2012}. Some double-sided inequalities of Ostrowski's type and some applications are also investigated 
(cf. \cite{MasDra14} and \cite{WAGVM}).  

In this paper, we present some selected Ostrowski type inequalities, considered from the point of view of numerical integration, precisely for error estimates in quadrature formulas. In Section 2 we give some basic facts on quadrature formulas of algebraic degree of exactness and approach for estimating the error term using Ostrowski type inequalities and Peano kernel techniques. In Section 3 we give several selected examples of simple quadrature formulas and the corresponding inequalities,  including the basic Ostrowski's inequality \cite{Os1938}, inequality of Milovanovi\'c and Pe\v cari\'c \cite{GM4}, Dragomir, Cerone and Roumeliotis \cite{DCR}, symmetric and asymmetric inequalities of Guessab and Schmeisser \cite{GSch} and Franji\'c 
\cite{IFr09}, respectively. Finally, we analyze  the four point symmetric inequality of Alomari \cite{Alom12}, as well as one variant with double internal nodes given by Liu and Park  \cite{LP}. 

\vspace{1.5cc}
\begin{center}
{\bf 2. PRELIMINARIES TO  QUADRATURE FORMULAS AND OSTROWSKI TYPE INEQUALITIES}
\end{center}

As we mentioned in Section 1, in 1938 Ostrowski proved the inequality (\ref{Os1}) 
for differentiable mappings with bounded first derivative. The constant ${1}/{4}$ in (\ref{Os1})  is sharp in the sense that it can not be replaced by a smaller one. 

In 1976, Milovanovi\'{c} and Pe\v{c}ari\'{c}  \cite{GM4}
presented the following generalization of Ostrowski's inequality with higher derivatives, i.e., when $\left\vert f^{\left( n\right)
}\left( x\right) \right\vert $ $\leq M$ $\left( \forall x\in \left(
a,b\right) \right) $ and $n>1$:
\medskip

\noindent
{\bf Theorem 2.1.} {\it
Let\ $f:\mathbb{R}\rightarrow \mathbb{R}$ \ be\ $n\left( >1\right) $
times differentiable function such that $\bigl\vert f^{(n)
}(x)\bigr\vert $ $\leq M$ $\bigl( \forall x\in (a,b) \bigr) $. 
Then, for every $x\in \left[ a,b\right] $%
\[\Biggl\vert \ \frac{1}{n}\Biggl( f(x) +\sum_{k=1}^{n-1}F_{k}\Biggr) -\frac{1}{b-a}\int_{a}^{b}f(t)
{\D}t\Biggr\vert \leq \frac{M}{n(n+1)!}\cdot \frac{(x-a)
^{n+1}+( b-x)^{n+1}}{b-a},\]
where $F_{k}$ is defined by 
\[
F_{k}\equiv F_{k}\left( f;n;x;a;b\right) \equiv \frac{n-k}{k!}\cdot \frac{
f^{( k-1) }( a) (x-a)^{k}-f^{(k-1)}(b)(x-b)^{k}}{b-a}.\]}

In a special case for $n=2$ and $|f''(x)|\le M$ on $(a,b)$, the previous inequality reduces to
\begin{align}\label{A-4}
\Biggl\vert \ \frac{1}{2}\Biggl(f(x) &+\frac{( x-a)
f( a) +( b-x) f( b) }{b-a}\Biggr) -\frac{1}{ 
b-a}\int_{a}^{b}f\left( t\right) {\D}t\,\Biggr\vert \\[2mm]
&\qquad\qquad\qquad\qquad\qquad\quad\leq \frac{M\left( b-a\right) ^{2}}{4}\left[ \frac{1}{12}+\frac{\left( x- 
\frac{a+b}{2}\right) ^{2}}{\left( b-a\right) ^{2}}\right] . \notag   
\end{align}

At the end of the nineties, there was an increased interest in this type of inequalities, and this increase has continued up to now. Many such integral inequalities for $n$-times differentiable mappings $(n\ge1)$  on the Lebesgue spaces $L^p[a,b], 1\le p\le +\infty$, have been obtained. Without loss of generality, we here consider some of these inequalities for functions given on $[-1,1]$, connected them to quadrature rules.  As usual the norm is defined by
\[{\|f\|}_p=\left\{\begin{array}{ll}
{\DS\left(\int_{-1}^1|f(t)|^p{\D}t\right)^{1/p}},&1\le p<+\infty,\\[5mm]
{\DS\sup_{t\in[-1,1]}|f(t)|},&p=+\infty,	
\end{array}\right.\]

In this way, the basic Ostrowski inequality (\ref{Os1}) becomes
\[\biggl\vert f(x)-\frac{1}{2}\int_{-1}^{1}f(t){\D}t\,\biggr\vert \leq  
 \frac{1+x^2}{4} \cdot2\left\Vert f'\right\Vert _{\infty },\quad x\in[-1,1],\]
i.e.,
 \begin{equation}\label{Osm1p1}
\biggl| \int_{-1}^{1}f(t){\D}t-2f(x)\biggr| \leq  (1+x^2)\left\Vert f'\right\Vert _{\infty },\quad x\in[-1,1].
\end{equation}
Similarly, the inequality (\ref{A-4}) reduces to 
 \begin{equation}\label{Os2m1p1}
\biggl| \int_{-1}^{1}f(t){\D}t-\frac{1}{2}\left[(1+x)f(-1)+2f(x)+(1-x)f(1)\right]\biggr|  \leq  \frac{1}{6}(1+3x^2)\left\Vert f''\right\Vert_{\infty },\end{equation}
for $x\in[-1,1]$. Now, these inequality (\ref{Osm1p1}) and (\ref{Os2m1p1}) can be treated  as  estimates of the remainder term of the one-point quadrature formula $Q_1(f)=2f(x)$ and the three-point quadrature formula 
 \begin{equation}\label{Os2m1p1QF}
Q_3(f)=\frac{1+x}{2}f(-1)+f(x)+\frac{1-x}{2}f(1), 
\end{equation}
respectively. In $Q_3(f)$ we have two fixed nodes $\pm1$ and one free $x$, $-1\le x\le 1$. In the case $x=\pm1$, $Q_3(f)$ reduces to the trapezoidal two-points formula $f(-1)+f(1)$.  

In general case we can consider $n$-point weighted  quadrature formulas 
\begin{equation}\label{qfQ}
I(wf)=\int_{-1}^1 w(t)f(t){\D}x=Q_n(f)+R_n(f),	
\end{equation} 
where $Q_n(f)$ is a quadrature sum, $R_n(f)$ is the corresponding remainder term, and $t\mapsto w(t)$ is a given weight function (for details see \cite[Sec.~5.1]{GMGVM08}).  Then  estimates of $|R_n(f)|$ lead to different  Ostrowski type inequalities in  certain classes of functions (for some collections of such  inequalities   see \cite{DragRas} and \cite{Drag17}). 

Let  $\PP_n$ be the set of all algebraic polynomials of degree at most $n$. 
The quadrature formula (\ref{qfQ}) has degree of exactness $d$ if for every $p\in \PP_d$ we have $R(p)=0$. In addition, if $R(p)\ne 0$ for some $p\in\PP_{d+1}$, this quadrature formula has precise degree of exactness $d$ (see Definition 5.1.2 in \cite[p.~320]{GMGVM08}). 

More generally, when a quadrature sum contains derivatives of arbitrary order at some points (nodes), such quadrature formulas are known as {\it Birkhoff-type quadratures} (cf. \cite{Shi2003}). The most important classes of such quadratures are ones with multiple nodes (for details see \cite{gvm2001,MSP2019} and the  references cited therein).

Here we consider only quadrature rules of the form  (Birkhoff type)
\begin{equation}\label{QF2}
Q_{n,m}(f)=\sum_{k=1}^n A_k f(x_k)+\sum_{k=1}^m B_k f'(y_k),
\end{equation}
for non-weighted integrals ($w(t)=1$), with the nodes $X_n=\{x_k\}_{k=1}^n$ and $Y_m=\{y_k\}_{k=1}^m$, such that
\[-1\le x_1<x_2<\cdots<x_n\le 1\quad \mbox{and}\quad -1\le y_1<y_2<\cdots<y_m\le 1.\]
These sets of nodes can have common points. In a special case it can be $n=m$ and $X_n=Y_m$, when we have a quadrature rule with multiple nodes. If the set $Y_m$  is empty, we have a standard quadrature formula with simple nodes.  

For a set of differentiable functions $\FF$, we define a linear functional $\LL:\FF\to\RR$ by means
\begin{equation}\label{LFun}
\LL f:=\int_{-1}^1 f(t){\D}t-\sum_{k=1}^n A_k f(x_k)-\sum_{k=1}^m B_k f'(y_k),
\end{equation}
and then we use the Peano representation of the functional $\LL f$, as well as a truncated power function  $t\mapsto (x-t)_+^r$, defined by
\[(x-t)_+^r=\left\{\begin{array}{ll}
(x-t)^r,& -\infty< t\le x,\\[2mm]
0,&t>x,	
\end{array}\right.\]
where $x$ is a fixed real number and $r$ is a nonnegative integer. Regarding (\ref{qfQ}) and (\ref{QF2}), we see that the remainder, in this case denoted by  $R_{n,m}(f)$, is  itself the  linear functional $\LL f$ on $\FF$. 

Suppose now that the quadrature formula (\ref{QF2}) has degree of exactness $d$ and that $f:[-1,1]\to\RR$ be a $(r+1)$-times differentiable function, where  $r\le d$. Then,  
according to the Peano kernel theorem \cite[Chp.~4]{GMP}, we have 
\begin{equation}\label{PKerTh}
R_{n,m}(f)=\int_{-1}^1 K_r(t) f^{(r+1)}(t){\D}t,
\end{equation}
where the $r$th Peano kernel is given by
\begin{equation}\label{PKer}
K_r(t)=\frac{1}{r!}	\LL(\,\cdot\,-t)_+^r.
\end{equation}
Applying H\"older's inequality to 
\[|R_{n,m}(f)|=\left|\int_{-1}^1 K_r(t) f^{(r+1)}(t){\D}t\right|,\]
with $1\le p\le+\infty$, $1/p+1/q=1$, we obtain 
\[\left|\int_{-1}^1 f(t){\D}t-Q_{n,m}(f)	\right|\le \bigl\|K_r\bigr\|_p \bigl\| f^{(r+1)}\bigr\|_q,\]
assuming that $f^{(r+1)}\in L^q[-1,1]$, where $q=p/(p-1)$.
\vspace{0.6cc}

In this paper we consider only the case when 
 $p=1$ $(q=+\infty)$, and it gives the following inequalities of Ostrowski type
\begin{equation}\label{PKerTh}
\left|\int_{-1}^1 f(t){\D}t-\sum_{k=1}^n A_k f(x_k)-\sum_{k=1}^m B_k f'(y_k)	\right|\le \left(\int_{-1}^1\left|K_r(t)\right|{\D}t\right) \bigl\| f^{(r+1)}\bigr\|_{\infty} 
\end{equation}
for  $r\le d$. Here, we need to determine $K_r(t)$ for the functional (\ref{LFun}), i.e.,
\begin{equation}\label{Kernel}  
r!K_r(t)=\frac{(1-t)^{r+1}}{r+1}-\sum_{x_k>t}A_k(x_k-t)^r-r\sum_{y_k>t}B_k(y_k-t)^{r-1},	
\end{equation}
as well as the integral of $\left|K_r(t)\right|$ over $[-1,1]$.

In the next section we analyze some typical inequalities of Ostrowski type (\ref{PKerTh}), starting with the basic inequality (\ref{Os1}), i.e., (\ref{Osm1p1}). To find the degree of exactness of a quadrature formula we check the values of the remainder term \begin{equation}\label{testR}
R_{n,m}(f)=I(f)-Q_{n,m}(f),
\end{equation}
taking the monomials $f(t)=e_k(t)=t^k$, $k=0,1,\ldots$\ .

\vspace{5.cc}
\begin{center}
{\bf 3. ANALYSIS OF CERTAIN TYPICAL INEQUALITIES OF OSTROWSKI TYPE}
\end{center}

\vspace{1.cc}

\noindent {\large\bf 3.1. Inequality of Ostrowski   \cite{Os1938}}
\vspace{0.5cc}

In this simplest case, for $n=1$, $m=0$, $A_1=2$, and $x_1=x$,  the quadrature formula (\ref{QF2}) reduces to $Q_{1,0}(f)=2f(x)$ with degree of exactness $d=0$, because of  $R_{1,0}(e_0)=I(e_0)-Q_{1,0}(e_0)=0$. Therefore,  the corresponding error estimate (\ref{PKerTh}) 
gives
\begin{equation}\label{Ost0}
\left|\int_{-1}^1f(t){\D}t-2f(x)\right|\le (1+x^2)\|f'\|_\infty\quad (-1\le x\le 1),	
\end{equation}
which is exactly Ostrowski's result (\ref{Osm1p1}). Here
\[K_0(t)=\left\{\begin{array}{ll}
\!-1-t,&-1<t\le x,\\[2mm]
\ \ 1-t,&\ \ x<t\le 1, 	
\end{array}\right.\quad\mbox{and}\quad \int_{-1}^1
\left|K_0(t)\right|{\D}t=1+x^2.\]

\begin{figure}[h]
\centering
\includegraphics[width=0.48\textwidth]{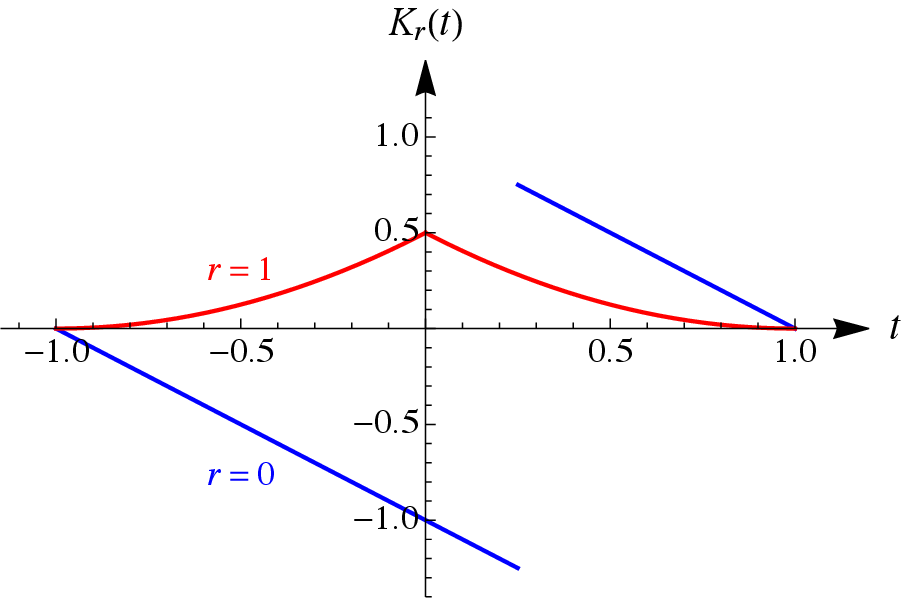}\ \ %
\includegraphics[width=0.48\textwidth]{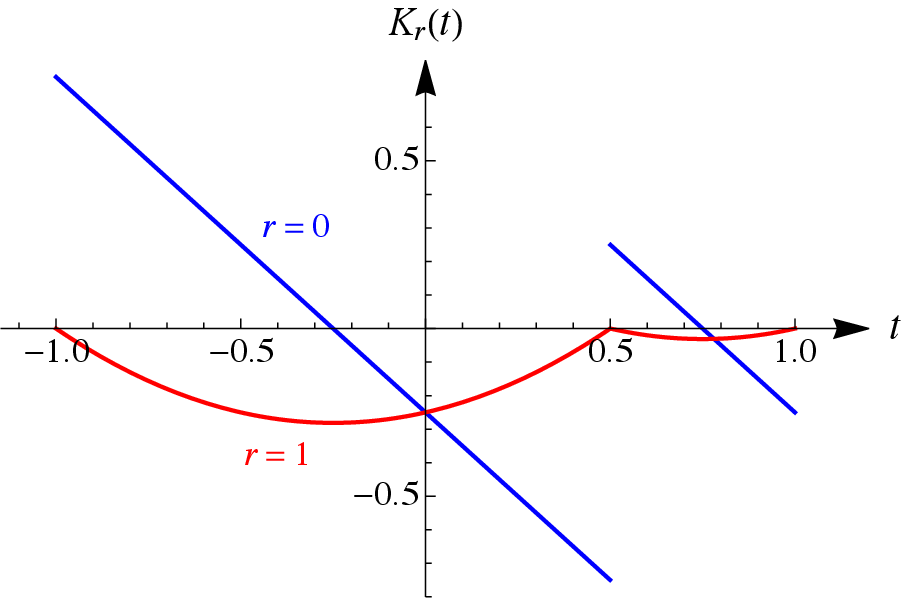}
\caption[]{The kernels $K_0(t)$ $(x=1/4)$ and $K_1(t)\ (x=0)$ in Ostrowski's inequality (left);   $K_0(t)$ and $K_1(t)$  for  $x=1/2$ in inequality of Milovanovi\'c and Pe\v cari\'c  (right)}\label{slOsMP01}	
\end{figure}

For $f(t)=e_1(t)$ and  $f(t)=e_2(t)$ we have $R_{1,0}(e_1)=I(e_1)-Q_{1,0}(e_1)=-2x$ and $R_{1,0}(e_2)=I(e_2)-Q_{1,0}(e_2)=\frac{2}{3}-2x^2$, for $x=0$ its degree of exactness becomes $d=1$, so that we have the following inequality
 \[\left|\int_{-1}^1f(t){\D}t-2f(0)\right|\le \frac{1}{3}\|f''\|_\infty ,\]
because
\[K_1(t)=\left\{\begin{array}{ll}
\dfrac{1}2(1+t)^2,&-1<t\le 0,\\[3mm]
\dfrac{1}2(1-t)^2,&\ \ 0<t\le 1, 	
\end{array}\right.\quad\mbox{and}\quad \int_{-1}^1
\left|K_1(t)\right|{\D}t=\frac{1}{3}.\]
In fact, this is an estimate for the one-point Gauss-Legendre formula
(cf. \cite[p.~172]{NAII}).

The kernel $t\mapsto K_0(t)$, for $x=1/4$, is presented in Figure \ref{slOsMP01} (left), as well as the kernel $t\mapsto K_1(t)$.
\vspace{.5cc}

\noindent {\large\bf 3.2. Inequality of Milovanovi\'c and Pe\v cari\'c \cite{GM4} }
\vspace{0.5cc}

The corresponding quadrature formula in (\ref{Os2m1p1}) is given by
(\ref{Os2m1p1QF}), where $n=3$ and $m=0$. Its nodes and weight coefficients are $-1$, $x$,  $1$  and $(1+x)/2$, $1$, ($1-x)/2$, respectively. According to $R_{3,0}(e_r)=I(e_r)-Q_{3,0}(e_r)=0$ for $r=0,1$, and $R_{3,0}(e_2)=-\frac{1}{3}-x^2\ne 0$, we see that degree of exactness is $d=1$.

Using (\ref{PKer}), i.e.,  (\ref{Kernel}), we determine the corresponding Peano kernels  for $r=0$ and $r=1$. Thus,  
\[K_0(t)=\left\{\begin{array}{ll}
\dfrac{1}2(x-2t-1),&-1<t\le x,\\[3mm]
\dfrac{1}2(x-2t+1),&\ \ x<t\le 1, 	
\end{array}\right.\!\!\qquad\]
and
\[K_1(t)=\left\{\begin{array}{ll}
\dfrac{1}2(1+t)(t-x),&-1<t\le x,\\[3mm]
\dfrac{1}2(-1+t)(t-x),&\ \ x<t\le 1, 	
\end{array}\right.\] 
as well as
\[\int_{-1}^1
\left|K_0(t)\right|{\D}t=\frac{1}{2}(1+x^2)   \quad\mbox{and}\quad    \int_{-1}^1
\left|K_1(t)\right|{\D}t=\frac{1}{6}(1+3x^2).\]
The kernels $t\mapsto K_0(t)$ and $t\mapsto K_1(t)$ for $x=1/2$ are presented in Figure \ref{slOsMP01} (right).

For $-1< x< 1$ we have the following inequalities of  Ostrowski's type
\[\left|\int_{-1}^1f(t){\D}t-f(x)-\frac{1}{2}\left[(1+x)f(-1)+(1-x)f(1)\right]\right|\le \frac{1}{2}(1+x^2)\|f'\|_\infty\]
and
\[\left|\int_{-1}^1f(t){\D}t-f(x)-\frac{1}{2}\left[(1+x)f(-1)+(1-x)f(1)\right]\right|\le \frac{1}{6}(1+3x^2)\|f''\|_\infty.\]
The second one is the original Milovanovi\'c-Pe\v cari\'c inequality (\ref{Os2m1p1}).
For $x=0$ the previous inequalities reduce to
\[\left|\int_{-1}^1f(t){\D}t-f(0)-\frac{1}{2}\left[f(-1)+f(1)\right]\right|\le \left\{\begin{array}{ll}
\dfrac{1}{2}	\|f'\|_\infty, & \mbox{when } r=0,\\[3mm]
\dfrac{1}{6}	\|f''\|_\infty, & \mbox{when } r=1.
\end{array}\right.\]
Otherwise, the last quadrature formula (for $x=0$) is a composition of two trapezoidal formulas,
\[Q_{3,0}(f)=f(0)-\frac{1}{2}\left[f(-1)+f(1)\right]=\frac{1}{2}\left[f(-1)+f(0)\right]+\frac{1}{2}\left[f(0)+f(1)\right].\]

Some generalizations of this kind of inequalities were given in \cite{DMP2000,DPU2003,PV2006}.
\vspace{.5cc}

\noindent {\large\bf 3.3. A  modification of the inequality (\ref{Os2m1p1})}
\vspace{0.5cc}

Now, we give a modification of (\ref{Os2m1p1}) by introducing a parameter $\lambda>0$.  Namely, instead of (\ref{Os2m1p1QF}), we 
consider a three-point quadrature rule of the form  
\begin{equation}\label{3rule}
Q_{3,0}^\lambda(f)=[1-\lambda(1-x)]f(-1)+2\lambda f(x)+[1-\lambda(1+x)]f(1),
\end{equation}
with  $-1 < x<1$, for which we have that
\[\left\{R_{3,0}^\lambda(e_r)\right\}_{r=0}^\infty=
\left\{0,0,-2 \lambda  \left(x^2-1\right)-\frac{4}{3},-2 \lambda  x \left(x^2-1\right),-2
   \lambda  \left(x^4-1\right)-\frac{8}{5},\ldots\right\}.\]
Note, that for $\lambda=0$, (\ref{3rule}) reduces to the well-known two-point trapezoidal rule.
   
The three-point quadrature rule (\ref{3rule}) has degree of exactness $d=1$, but for $\lambda=\lambda_x=\frac{2}{3}(1-x^2)^{-1}$  and  $x\ne 0$, the rule 
(\ref{3rule}) reduces to
\begin{equation}\label{3rule2}
Q_{3,0}^{\lambda_x}(f)=\frac{1}{3}\left[\frac{1+3x}{1+x}f(-1)+\frac{4}{1-x^2}f(x)+\frac{1-3x}{1-x}f(1)\right],  
\end{equation}
with degree of exactness $d=2$, because $R(e_r)=0$, $r=0,1,2$, and $R(e_3)=4x/3$. This kind of quadrature  rules have been also treated in \cite[\S6.2]{Mon2012}.
If we need to have a quadrature rule with all positive weight coefficients, then the parameter $x$ must be 
$|x|<1/3$.
\begin{figure}[h]
\centering
\includegraphics[width=0.48\textwidth]{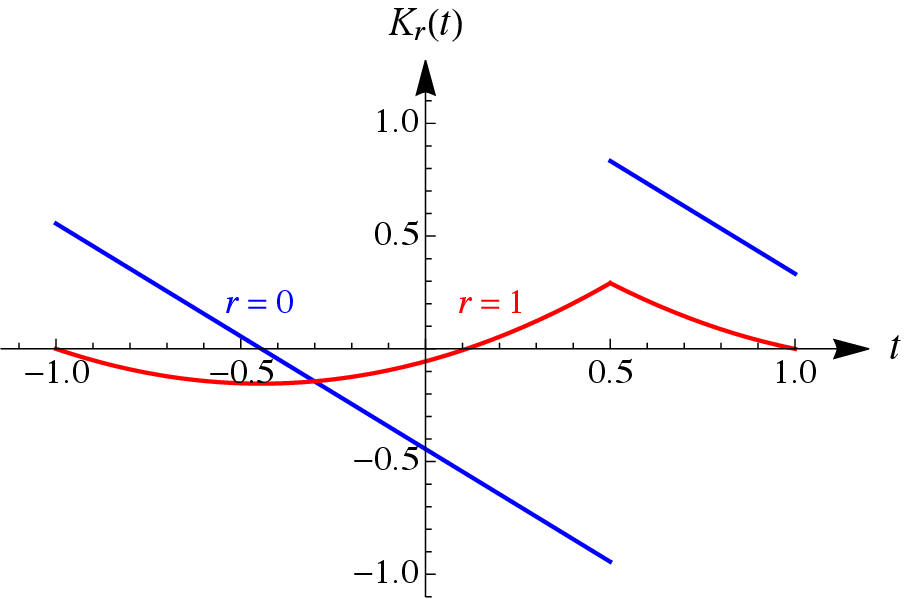}\ %
\includegraphics[width=0.48\textwidth]{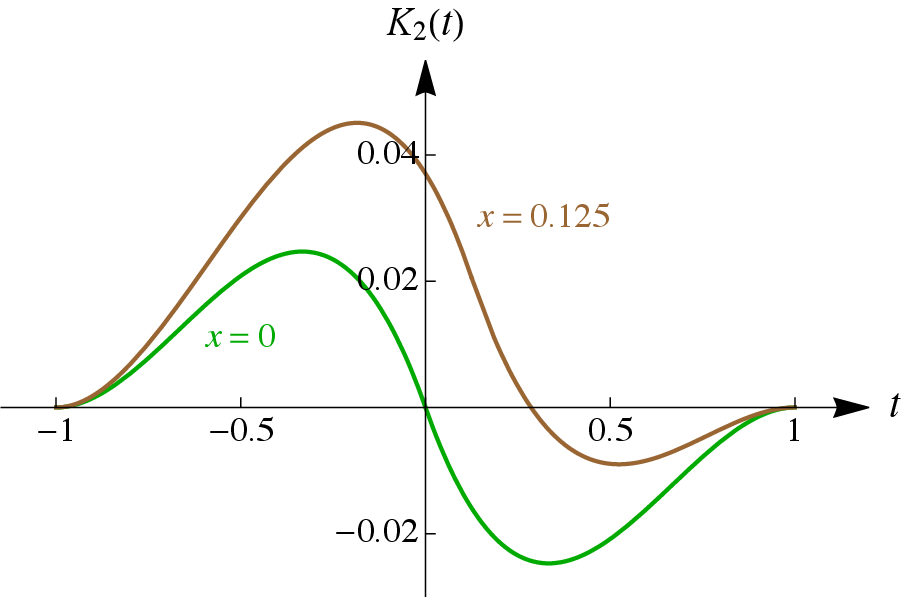}
\caption[]{The  kernels  $K_0(t)$ and $K_1(t)$ for $x=1/2$ (left) and $K_2(t)$ for $x=0$ and $x=0.125$ (right)}\label{slgvm0123}	
\end{figure}

The corresponding kernels $K_r(t)$ for $r=0,1,2$ are given by
\[K_0(t)=\left\{\begin{array}{ll}
-t-\dfrac{2}{3 (1+x)},&-1<t\le x,\qquad\qquad\quad\\[4mm]
-t+\dfrac{2}{3 (1-x)},&\ \ x<t\le 1, 	
\end{array}\right. \]

\[K_1(t)=\left\{\begin{array}{ll}
\dfrac{(1+t) (3 t x+3 t-3 x+1)}{6 (1+x)},&-1<t\le x,\\[3mm]
\dfrac{(1-t) (3 t x-3 t+3 x+1)}{6 (1-x)},&\ \ x<t\le 1, 	
\end{array}\right.\]
\[K_2(t)=\left\{\begin{array}{ll}
\dfrac{(1+t)^2 [2x-t (1+x)]}{6(1+x)},&-1<t\le x,\\[4mm]
\dfrac{(1-t)^2 [2x-t (1-x)]}{6(1-x)},&\ \ x<t\le 1, 	
\end{array}\right.\]
respectively. These kernels are displayed in Fig.~\ref{slgvm0123}.

\begin{figure}[h]
\centering
\includegraphics[width=0.5\textwidth]{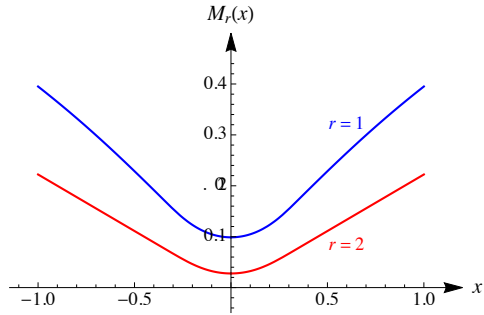} 
\caption[]{The bound function $x\mapsto M_r(x)$ for $r=1$ and $r=2$}\label{slMr12}	
\end{figure}

For $0\le x<1$, we have 
\[M_0(x)=\int_{-1}^1
\left|K_0(t)\right|{\D}t=\left\{\begin{array}{ll}
\dfrac{9 x^6+3 x^4-x^2+5}{9 \left(1-x^2\right)^2},&0\le x<\dfrac{1}{3},\ \quad\\[4mm]
\dfrac{\left(3 x^2+3 x+2\right)^2}{9 (1+x)^2},&\dfrac{1}{3}\le x<1,
\end{array}\right.\]

\[M_1(x)=\int_{-1}^1
\left|K_1(t)\right|{\D}t=\left\{\begin{array}{ll}
\dfrac{8 \left(1-3 x^2\right) \left(3 x^2+1\right)^2}{81 \left(1-x^2\right)^3},&0\le x<\dfrac{1}{3},\\[4mm]
\dfrac{4 (3 x+1)^3}{81 (1+x)^3},&\dfrac{1}{3}\le x<1,
\end{array}\right.\]

\[M_2(x)=\int_{-1}^1
\left|K_2(t)\right|{\D}t=\left\{\begin{array}{ll}
\dfrac{8 x^5+49 x^4-60 x^3+22 x^2-4 x+1}{36 (1-x)^4},&0\le x<\dfrac{1}{3},\\[4mm]
\dfrac{2x}{9},&\dfrac{1}{3}\le x<1,
\end{array}\right.\]
and $M_k(-x)=M_k(x)$. 
 These bounds $x\mapsto M_r(x)$ for $r=1$ and $r=2$ are presented in Fig.~\ref{slMr12}. 
 
For example, for $r=2$ we have the following inequality
\[\left|\int_{-1}^1f(t){\D}t-\frac{1}{3}\left[\frac{1+3x}{1+x}f(-1)+\frac{4}{1-x^2}f(x)+\frac{1-3x}{1-x}f(1)\right]\right|\le M_2(x)\|f'''\|_\infty.\]

In the case  $x=0$,  the quadrature formula (\ref{3rule2})  becomes the well-known Simpson rule
\begin{equation}\label{3rule3}
Q(f)=\frac{1}{3}\left[f(-1)+4f(0)+f(1)\right],  
\end{equation}
with degree  of exactness $d=3$ and its kernel is
\begin{equation}\label{SimpK3}
K_3(t)=\left\{\begin{array}{ll}
\ \ \dfrac{1}{72}(1+t)^3(3t-1),&-1<t\le 0,\\[4mm]
-\dfrac{1}{72}(1-t)^3(3t+1),&\ \ 0<t\le 1, 	
\end{array}\right.	
\end{equation} 
with  the bound constant
\begin{equation}\label{SimpK3M3} 
M_3=\int_{-1}^1
\left|K_3(t)\right|{\D}t	=\frac{1}{90}.
\end{equation}

\noindent {\large\bf 3.4. Inequality of Dragomir, Cerone and Roumeliotis \cite{DCR} }
\vspace{0.5cc}

A similar formula to (\ref{3rule}) was considered   
by Dragomir, Cerone and Roumeliotis \cite{DCR} in the form
\begin{equation}\label{DCR2000}
\widetilde Q_{3,0}^\lambda(f)= \lambda[f(-1)+f(1)]+2(1-\lambda) f(x),	
\end{equation}
for all $\lambda\in[0,1]$ and $-1+3\lambda/2\le x\le 1-3\lambda/2$. According to (\ref{testR}) we have
\begin{align*}
\left\{\widetilde R_{3,0}^\lambda(e_r)\right\}_{r=0}^\infty=&
\left\{0,\,2(\lambda-1)x,\,-2\lambda+2(\lambda -1)x^2+\frac{2}{3},\,2 (\lambda -1) x^3,\right.\\
& \left.\  -2 \lambda +2 (\lambda -1)
   x^4+\frac{2}{5},\, 2 (\lambda -1) x^5,\,-2 \lambda +2 (\lambda -1) x^6+\frac{2}{7},\,  \ldots\right\}.	
\end{align*}
If $\lambda\ne1$ and $x\ne0$,  we see that the rule (\ref{DCR2000}) has degree of exactness $d=0$ and we have  
\begin{equation}\label{SimpK0}
 K_0(t)=\left\{\begin{array}{ll}
\!-1-t+\lambda,&-1<t\le x,\\[2mm]
\ \ 1-t-\lambda,&\ \ x<t\le 1, 	
\end{array}\right.
\end{equation} 
and
\[
M_0(x)= \int_{-1}^1
\left|K_0(t)\right|{\D}t=\lambda^2+(1-\lambda)^2+x^2.\]
Dragomir, Cerone and Roumeliotis \cite{DCR}
 obtained the following inequality
\begin{equation}\label{DCR2000Ineq}
\left|\int_{-1}^1f(t){\D}t-\widetilde Q_{3,0}^\lambda(f)\right|\le \left[\lambda^2+(1-\lambda)^2+x^2\right]\|f'\|_\infty \end{equation}
for differentiable functions $f:[-1,1]\to \RR$ with bounded derivative on $(-1,1)$.

Some similar inequalities were obtained by Ujevi\'c \cite{UJ2003}.

Evidently, for $\lambda=0$ (one-point rule) the inequality reduces to  Ostrowski's inequality (\ref{Ost0}), while for $\lambda=1$ (two-point trapezoidal rule, with $d=1$) it gives 
\[\left|\int_{-1}^1f(t){\D}t-f(-1)-f(1)\right|\le \|f'\|_\infty.\]

For $\lambda=1/3$ and $x\in[-1/4,1/4]$, (\ref{DCR2000Ineq}) reduces to the generalized Simpson inequality
\begin{equation}\label{Simpx}
\left|\int_{-1}^1f(t){\D}t-\frac{1}{3}\left[f(-1)+4f(x)+f(1)\right]\right|\le \left(\frac{5}{9}+x^2\right)\|f'\|_\infty. 	\end{equation}

For $\lambda=1/3$ and  $x=0$, the quadrature rule (\ref{DCR2000}) becomes the standard Simpson formula (\ref{3rule3}), with degree of exactness  now $d=3$, and we can give the error estimates for each $r\le d=3$.  

The kernels $K_r(t)$ for $r=1$ and $r=2$ are
\[K_1(t)=\left\{\begin{array}{ll}
\dfrac{1}{6}(1+4t+3t^2),&-1\le t\le 0,\\[4mm]
\dfrac{1}{6}(1-4t+3t^2),&\ \ \,0<t\le 1,	
\end{array}\right.\ 
K_2(t)=\left\{\begin{array}{ll}
-\dfrac{1}{6}t(1+t)^2,&-1\le t\le 0,\\[4mm]
-\dfrac{1}{6}t(1 - t)^2,&\ \ \,0<t\le 1,	
\end{array}\right.\] 
so that the corresponding  bounds are
\begin{equation}\label{grM1M2}
M_1=\int_{-1}^1
\left|K_1(t)\right|{\D}t=\frac{8}{81}\quad \mbox{and}\quad 
  M_2=\int_{-1}^1
\left|K_2(t)\right|{\D}t=\frac{1}{36}.	
\end{equation}

\begin{figure}[h]
\centering
\includegraphics[width=0.48\textwidth]{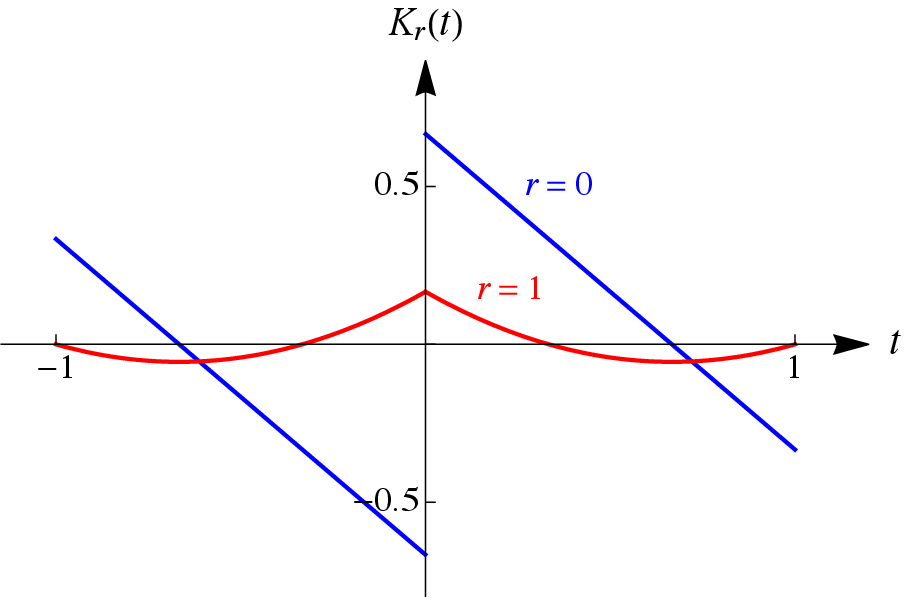}\quad
\includegraphics[width=0.48\textwidth]{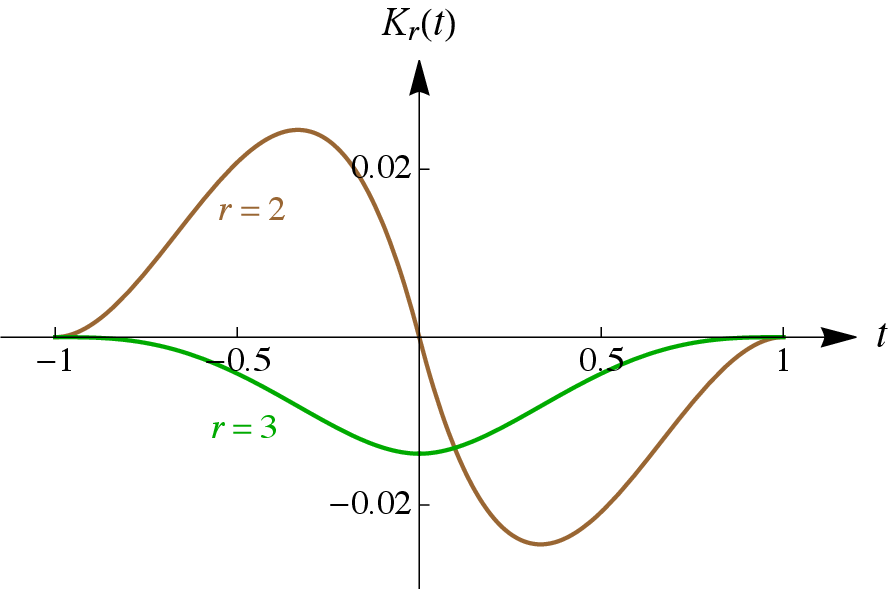}
\caption[]{The  kernels of the Simpson rule $K_r(t)$: (left) $r=0$ and $r=1$;  (right) $r=2$ and $r=3$}\label{SimpKer}.	
\end{figure}

The  kernels $K_0(t)$ (given by (\ref{SimpK0}) for $\lambda=1/2$ and $x=0$), $K_1(t)$, $K_2(t)$ and $K_3(t)$ (given earlier by  
(\ref{SimpK3})) are displayed  in Fig.~\ref{SimpKer}.

Thus, for the Simpson formula, according to (\ref{Simpx}) for $x=0$, previous bounds (\ref{grM1M2}), as well as (\ref{SimpK3M3}), we have the following estimates    
\[
\left|\int_{-1}^1f(t){\D}t-\frac{1}{3}\left[f(-1)+4f(0)+f(1)\right]\right|\le 
\left\{\begin{array}{ll}
\dfrac{5}{9}\|f'\|_\infty, & \mbox{when } r=0,\\[3mm]
\dfrac{8}{81}\|f''\|_\infty, & \mbox{when } r=1,\\[3mm]
\dfrac{1}{36}\|f'''\|_\infty, & \mbox{when } r=2,\\[3mm]
\dfrac{1}{90}\|f^{iv}\|_\infty, & \mbox{when } r=3.
\end{array}\right.\]

\noindent {\large\bf 3.5. Symmetric inequality of Guessab and Schmeisser \cite{GSch}}
\vspace{0.5cc}

Here we have a symmetric quadrature rule $Q_{2,0}(f)=f(-x)+f(x)$
 (see \cite{GSch}). Suppose that $0\le x\le 1$. According to (\ref{testR}) we get
 \[\left\{R_{2,0}^\lambda(e_r)\right\}_{r=0}^\infty=
\left\{0,\,0,\,\frac{2}{3}-2 x^2,\,0,\,\frac{2}{5}-2 x^4,\,0,\,\ldots\,\right\},\]
from which we conclude that the rule $Q_{2,0}(f)$ has degree of exactness $d=1$ for each $x\ne 1/\sqrt{3}$. For $x=1$ this rule reduces to the trapezoidal rule, given also as a special case of
$\widetilde Q_{3,0}^\lambda(f)$  (Eq. (\ref{DCR2000}) in \$3.4) 
for $\lambda=1$.  

However,  for  $x=1/\sqrt{3}$ this degree of exactness becomes $d=3$. 

It is easy to find the kernels for $r=0$ and $r=1$,
\[K_0(t)=\left\{\begin{array}{ll}
-1-t,&-1\le t\le -x,\\[2mm]
-t,&-x<t\le x,\\[2mm]
1-t,&\ \ \,x<t\le 1,	
\end{array}\right.\quad 
K_1(t)=\left\{\begin{array}{ll}
\frac{1}{2}(1+t)^2,&-1\le t\le -x,\\[2mm]
\frac{1}{2}(1 + t^2 - 2 x),&-x<t\le x,\\[2mm]
\frac{1}{2}(1 - t)^2,&\ \ \,x<t\le 1,	
\end{array}\right.\]
for which the bounds are even functions,  given by 
\[M_0(x)=\int_{-1}^1
\left|K_0(t)\right|{\D}t=1-2x+2x^2\quad (0\le x\le 1)\]
and 
\[M_1(x)=\int_{-1}^1
\left|K_1(t)\right|{\D}t=\left\{\begin{array}{ll}
\dfrac{1}{3}-x^2,&0\le x \le\dfrac{1}{2},\\[4mm]
\dfrac{1}{3}\left[ 4 (2 x-1)^{3/2}+1-3 x^2\right],&\dfrac{1}{2}<x\le 1.
\end{array}\right.\]
\begin{figure}[h]
\centering
\includegraphics[width=0.48\textwidth]{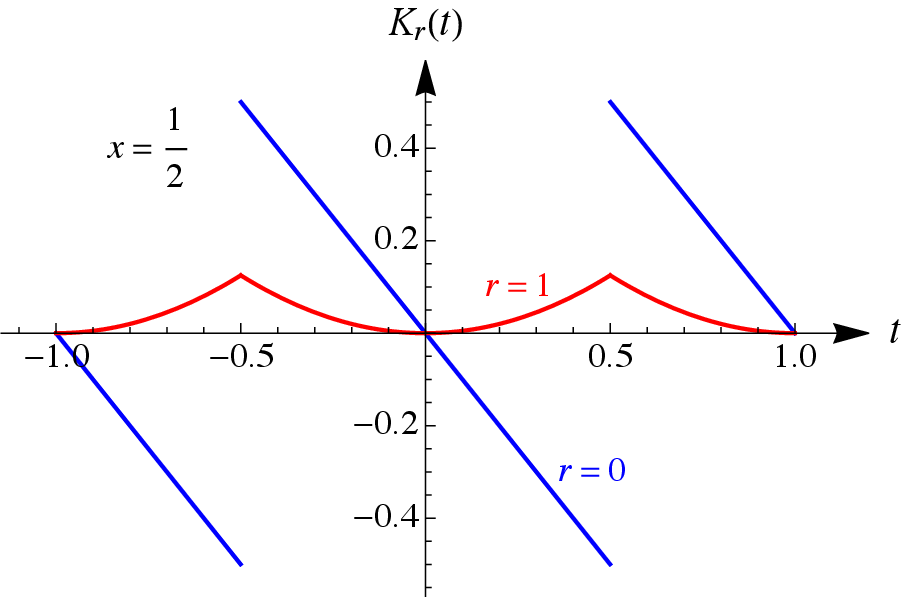}\quad
\includegraphics[width=0.48\textwidth]{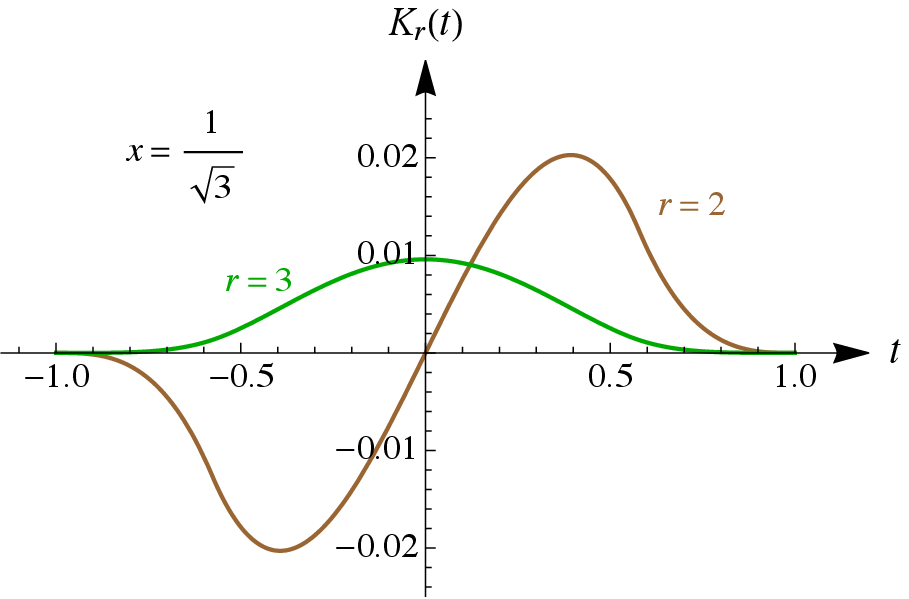}
\caption[]{The  kernels $K_r(t)$: (left) $K_0(t)$ and $K_1(t)$,  when  $x=1/2$;  (right)  $K_2(t)$ and $K_3(t)$ for $x=1/\sqrt{3}$}\label{slAlGe0123}	
\end{figure}

Thus,
\[\left|\int_{-1}^1f(t){\D}t-\left[f(-x)+f(x)\right]\right|\le 
M_r(x)\|f^{(r+1)}\|_\infty\quad (r=0,1),\]
where $M_0(x)$ and $M_1(x)$ are given above. It is interesting to mention that
\[\min_{x\in(0,1]}M_0(x)=M_0\Bigl(\frac{1}{2}\Bigr)=\frac{1}{2}\quad\mbox{and}\quad
\min_{x\in(0,1]}M_1(x)=M_1\bigl(4-2\sqrt{3}\bigr)=7-4\sqrt{3}.\]
The last optimal value was also obtained in \cite{UJ2004}.

For $x=1/\sqrt{3}$, quadrature formula $Q_{2,0}(f)$ \cite{GSch} reduces to the two-point Gauss-Legendre formula with degree of precision $d=3$  (cf. \cite[\S7.2.9]{NAII}). In that case for the kernels $K_2(t)$ and $K_3(t)$ we have
\[K_2(t)=\left\{\begin{array}{ll}
-\dfrac{1}{6}(1+t)^3,&-1\le t\le -\dfrac{1}{\sqrt3},\\[4mm]
-\dfrac{t}{6}\left(3-2\sqrt{3}+t^2\right),&-\dfrac{1}{\sqrt3}<t\le \dfrac{1}{\sqrt3},\qquad\qquad\qquad\qquad\quad\\[4mm]
\ \ \,\dfrac{1}{6}(1-t)^3,&\ \ \,\dfrac{1}{\sqrt3}<t\le 1,	
\end{array}\right.\]
and
\[K_3(t)=\left\{\begin{array}{ll}
\dfrac{1}{24}(1+t)^4&-1\le t\le -\dfrac{1}{\sqrt3},\\[4mm]
\dfrac{1}{216}\left(9 t^4-18 \left(2 \sqrt{3}-3\right) t^2-4 \sqrt{3}+9\right),&-\dfrac{1}{\sqrt3}<t\le \dfrac{1}{\sqrt3},\\[4mm]
\dfrac{1}{24}(1 - t)^4,&\ \ \,\dfrac{1}{\sqrt3}<t\le 1,	
\end{array}\right.\]
respectively, with the bounds in the corresponding two-point Gauss-Legendre rule $M_2\|f'''\|_\infty$ and $M_3\|f^{iv}\|_\infty$, where 
\[M_2=\int_{-1}^1 \left|K_2(t)\right|{\D}t=\frac{1}{108}\left(9-4\sqrt{3}\right)\quad\mbox{and}\quad
M_3=\int_{-1}^1 \left|K_3(t)\right|{\D}t=\frac{1}{135}.\]

A general two-point integral quadrature formula, using the concept of harmonic polynomials, was derived in \cite{KPVNonlAnal08}.

\vspace{0.8cc}

\noindent {\large\bf 3.6. Asymmetric inequality of Franji\'c \cite{IFr09}}
\vspace{0.5cc}

Instead of symmetric rule $Q_{2,0}(f)=f(-x)+f(x)$, Frani\'c \cite{IFr09} considered asymmetric rules with a fixed node (Radau type), using the extended Euler formula  obtained earlier in \cite{DMP2000}. 

\begin{figure}[h]
\centering
\includegraphics[width=0.48\textwidth]{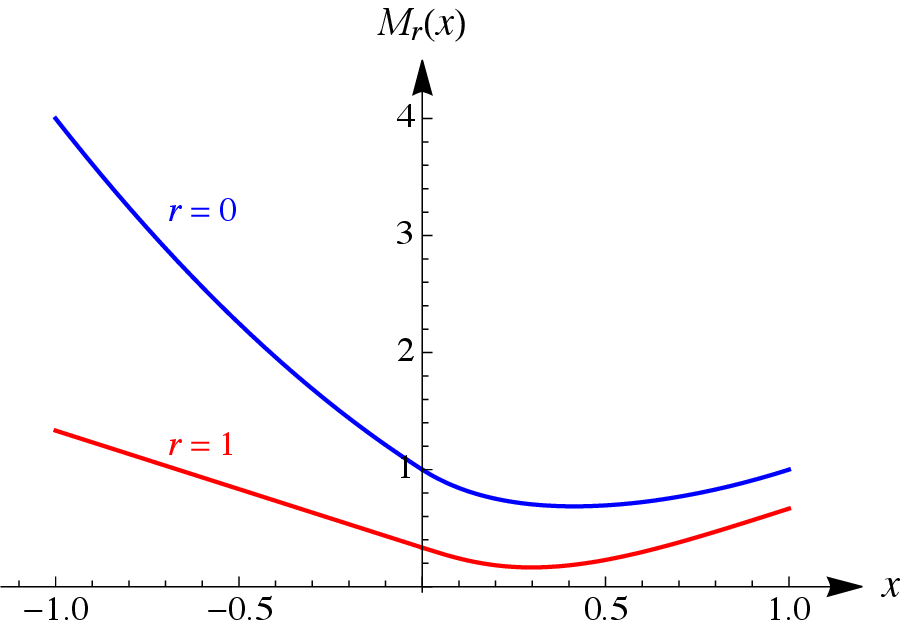} 
\includegraphics[width=0.48\textwidth]{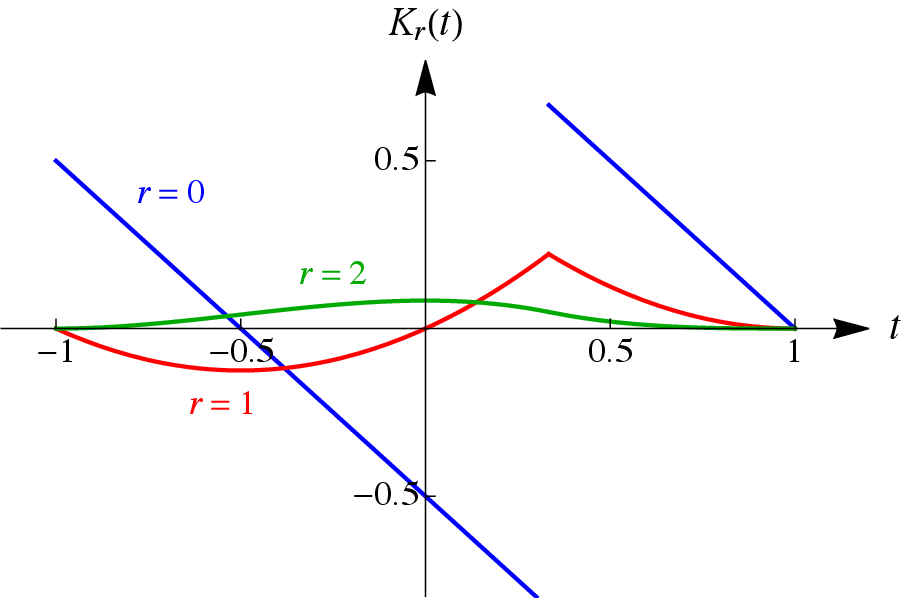} 
\caption[]{The bound function $x\mapsto M_r(x)$ for $r=0$ and $r=1$ (left) and the  kernels $K_r(t)$, $r=0,1,2$ for the Radau rule (\ref{Radau1})}\label{Rad1Mr01}	
\end{figure}

Here we fix the end-point $-1$ and consider the rule
\[Q_{2,0}^\lambda(f)=\lambda f(-1)+(2-\lambda)f(x),\]
with positive weight coefficients, i.e.,  when $0<\lambda<2$. According to (\ref{testR}) we have
\[\left\{R_{2,0}^\lambda(e_r)\right\}_{r=0}^\infty=
\left\{0,\,\lambda +(\lambda -2) x,\,-\lambda +(\lambda -2) x^2+\frac{2}{3},\,\lambda +(\lambda -2) x^3,\,\ldots\,\right\}.\]
Taking $\lambda=\lambda_x=2x/(1 + x)$, we have that
\[\left\{R_{2,0}^{\lambda_x}(e_r)\right\}_{r=0}^\infty=
\left\{0,\,0,\,2\left(\frac{1}{3}- x\right),-2 (x-1) x,\, \ldots\,\right\},\]
i.e., $Q_{2,0}^{\lambda_x}(f)$ is a rule of degree of the exactness $d=1$, except the case  $x\ne1/3$. The kernels of this  quadrature rule,
\begin{equation}\label{GenFix}
Q_{2,0}^{\lambda_x}(f)=\frac{2 x}{1 + x}f(-1)+\frac{2}{1 + x}f\left(x\right), 	
\end{equation} 
are
\begin{equation}\label{Rad0}
K_0(t)=\left\{\begin{array}{ll}
1-t-\dfrac{2}{1 + x},&-1\le t\le x,\\[4mm]
1 - t,&\ \ \,x<t\le 1,	
\end{array}\right.\qquad\qquad\quad	
\end{equation} 
and 
\begin{equation}\label{Rad1}
K_1(t)=\left\{\begin{array}{ll}
\dfrac{(t+1) (t x+t-3 x+1)}{2 (x+1)},&-1\le t\le x,\\[4mm]
-\dfrac{1}{2}t(1 - t)^2,&\ \ \,x<t\le 1,	
\end{array}\right.	
\end{equation} 
with the bounds $\DS M_k(x)=\int_{-1}^1|K_r(t)|{\D}t$, $r=0,1$, given by 
\[M_0(x)=\left\{\begin{array}{ll}
\!\!(1-x)^2,&\!\!\!\!\!-1<x\le 0,\\[3mm]
\!\!\dfrac{(1+x^2)^2}{(1+x)^2},&0<x\le 1,	
\end{array}\right.  
M_1(x)=\left\{\begin{array}{ll}
\!\dfrac{1}{3}(1-3x),&\!\!\!\!\!-1<x\le 0,\\[3mm]
\!\dfrac{1-6 x^2+24 x^3-3 x^4}{3 (x+1)^3},&0<x\le 1,	
\end{array}\right.\]
and presented in Fig.~\ref{Rad1Mr01} (left). Their minimal values  on $(0,1]$ are 
\[\min_{x\in(0,1]}M_0(x)=M_0\bigl(\sqrt{2}-1\bigr)=12 - 8\sqrt{2}\approx 0.686292\]
and
\begin{align*}
\min_{x\in(0,1]}M_1(x)=&M_1\Bigl(2 \sqrt{2}-1-2 \sqrt{2-\sqrt{2}}\,\Bigr)\\
=&\frac{4}{3} \left(5-3 \sqrt{2}-2 \sqrt{10-7 \sqrt{2}}\,\right)\approx 0.164412	
\end{align*}
(see also \cite[Theorem 3]{IFr09} and  \cite[pp.~253--254]{Mon2012}).

For $x=1/3$ this asymmetric rule (\ref{GenFix}) reduces to the simplest Radau formula
\begin{equation}\label{Radau1}
Q_{2,0}^{1/2}(f)=\frac{1}{2}f(-1)+\frac{3}{2}f\left(\frac{1}{3}\right),\end{equation} 
with degree of exactness $d=2$. Its kernels $K_0(t)$ and $K_1(t)$
are given by (\ref{Rad0}) and (\ref{Rad1}), respectively, for $x=1/3$. For the kernel $K_2(t)$ we obtain
\[K_2(0)= \left\{\begin{array}{ll}
\dfrac{1}{12} (1-2 t) (1+t)^2,&-1\le t\le \dfrac{1}{3},\\[4mm]
\dfrac{1}{6} (1-t)^3,&\ \ \,\dfrac{1}{3}<t\le 1.
\end{array}\right.	
\]
These kernels are displayed in Fig.~\ref{Rad1Mr01} (right), and the corresponding inequalities are
\[\left|\int_{-1}^1f(t){\D}t-\left[\frac{1}{2}f(-1)+\frac{3}{2}f\left(\frac{1}{3}\right)\right]\right|\le 
\left\{\begin{array}{ll}
\dfrac{25}{36}\|f'\|_\infty, & \mbox{when } r=0,\\[3mm]
\dfrac{1}{6}\|f''\|_\infty, & \mbox{when } r=1,\\[3mm]
\dfrac{2}{27}\|f'''\|_\infty, & \mbox{when } r=2.
\end{array}\right.\]

For some other inequalities of this type see \cite{IFr09,Mon2012}.

Now we mention also a general two-point rule
\[\widehat Q_{2,0}^\lambda(f)=(1+\lambda) f(x)+(1-\lambda)f(y),\]
with $-1\le x\le \lambda \le y\le 1$, considered by Alomari \cite{Alom17}. Since
 \begin{align*}
 \left\{\widehat R_{2,0}^\lambda(e_r)\right\}_{r=0}^\infty=&
\left\{0,\,(\lambda -1) y-(\lambda +1) x,\,-(\lambda +1) x^2+(\lambda -1) y^2+\frac{2}{3}\right.,\\
&\ \ \left.(\lambda -1) y^3-(\lambda +1) x^3,\,-(\lambda
   +1) x^4+(\lambda -1) y^4+\frac{2}{5},\ldots\right\},
 \end{align*}
we conclude that $d=0$, except the case when  $\lambda=\lambda_{x,y}=(x+y)/(y-x)$.

In this general case Alomari \cite{Alom17} obtained the following inequality
\begin{align*}
&\left|\int_{-1}^1f(t){\D}t-\left[(1+\lambda) f(x)+(1-\lambda)f(y)\right]\right|\\[2mm]
&\qquad\le \left\{\frac{1}{4}\left[(1-\lambda )^2+(1+\lambda)^2\right]+\left(x+\frac{1-\lambda }{2}\right)^2+\left(y-\frac{1+\lambda}{2}\right)^2\right\}\|f'\|_\infty, 	
\end{align*}
as well as several particular cases of this inequality.

\vspace{0.6cc}

\noindent {\large\bf 3.7. Four point symmetric inequality of Alomari \cite{Alom12}}
\vspace{0.5cc} 

Inspired by (\ref{DCR2000}), Alomari \cite{Alom12} considered 
the following symmetric rule
\begin{equation}\label{Alom12}
 Q_{4,0}^\lambda(f)= \lambda[f(-1)+f(1)]+(1-\lambda) [f(-x)+f(x)]\end{equation}
in three different classes:  (a)  functions of bounded variation, (b) absolutely continuous functions whose first derivative belongs 
to $L^\infty[-1,1]$, and (c) absolutely continuous functions whose first derivative belongs to $L^p[-1,1]$ for $p>1$. Following our previous discussion, we interested only in   rules and corresponding Ostrowski type inequalities for (sufficiently differentiable) functions with  $f^{(r+1)}\in L^\infty[-1,1]$, where $r\le d$ and $d$ is degree of exactness of the rule $Q_{4,0}^\lambda(f)$. We suppose that $\lambda\in(0,1)$ in order to have positive weight coefficients, as well as $0\le x\le 1$. Note that for $\lambda=0$ and $\lambda=1$ the rule (\ref{Alom12}) reduces to $Q_{2,0}(f)$ from \S3.5.

According to (\ref{testR}) we get
 \begin{align*}
 \left\{R_{4,0}^\lambda(e_r)\right\}_{r=0}^\infty=&
\left\{0,\,0,\,-2 \lambda +2 (\lambda -1) x^2+\frac{2}{3},\,0,\,-2 \lambda +2 (\lambda -1) x^4+\frac{2}{5},\,0,\right.\\ 	
&\left.\quad-2 \lambda +2
   (\lambda -1) x^6+\frac{2}{7},\,0,\,-2 \lambda +2 (\lambda -1) x^8+\frac{2}{9},\,0,\, \ldots\, \right\}
 \end{align*}
and conclude that  this rule has degree of exactness $d=1$, except the cases when $x^2=(1/3-\lambda)/(1-\lambda)\ne 1/5$, i.e.,   
\[\lambda=\lambda_x=\frac{1-3 x^2}{3 \left(1-x^2\right)}\ne\frac{1}{6}, \]
when $d=3$. Then, we have
\[\left\{R_{4,0}^{\lambda_x}(e_r)\right\}_{r=0}^\infty= 
\left\{0,\,0,\,0,\,0,\,\frac{4}{15} \left(5 x^2-1\right),\,0,\,\frac{4}{21} \left(7 x^4+7 x^2-2\right),\,0,\,\ldots\,\right\}.\] 
Finally, for $x=\xi=1/\sqrt{5}$ and $\lambda=1/6$, it reduces to 
\[\left\{R_{4,0}^{1/6}(e_r)\right\}_{r=0}^\infty= 
\left\{0,\,0,\,0,\,0,\,0,\,0,\,-\frac{32}{525},\,0,\,\ldots\,\right\}\quad\mbox{and}\quad d=5.\]

\begin{figure}[h]
\centering
\includegraphics[width=0.48\textwidth]{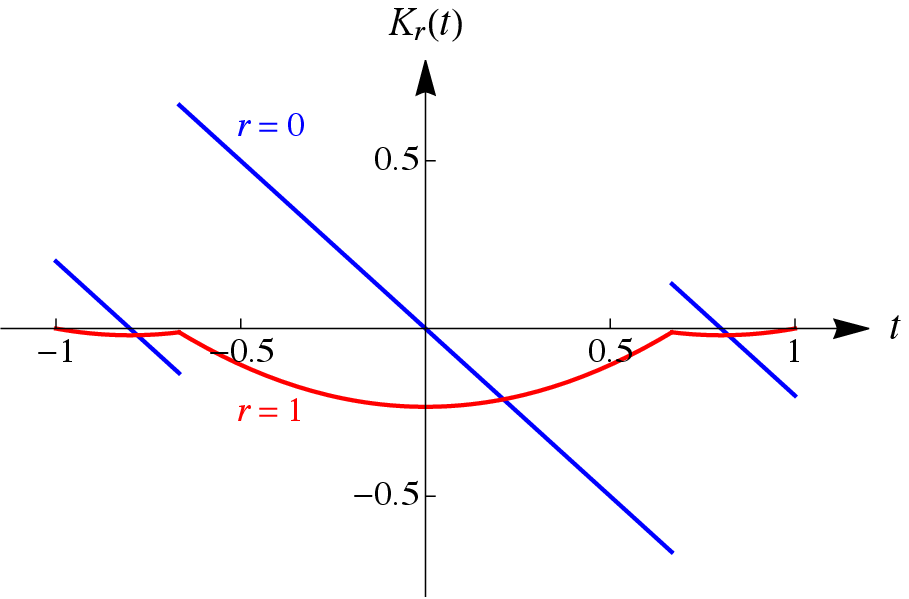}\quad
\includegraphics[width=0.48\textwidth]{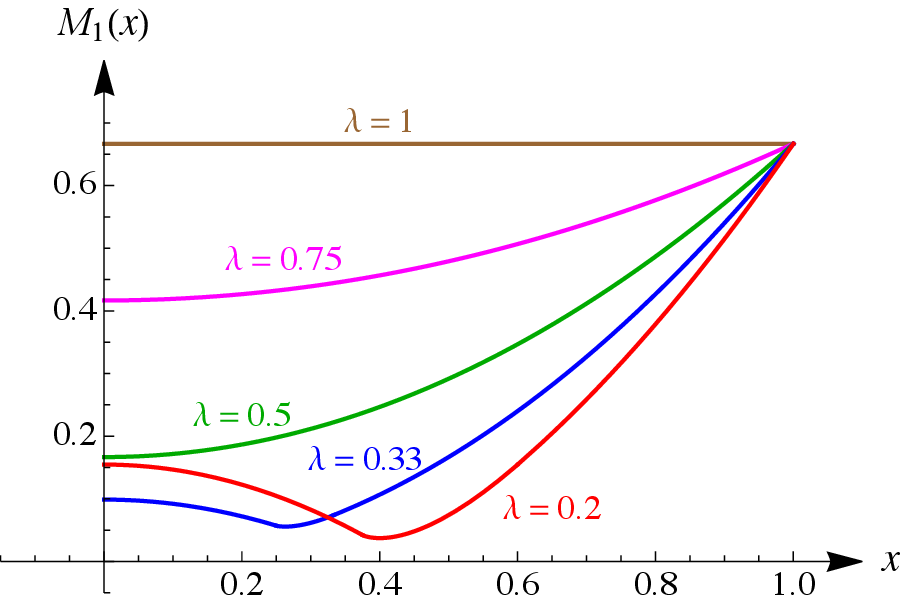}
\caption[]{The  kernels $K_r(t)$ for (\ref{Alom12}): (left) $K_0(t)$ and $K_1(t)$,  when $\lambda=1/5$ and  $x=2/3$;  (right)  The bound $M_1(x)$ for different values of $\lambda$}\label{AlomFig1}	
\end{figure}

In the case when $d=1$,  the kernels $K_r(t)$ $(r=0,1)$ are given by
\begin{equation}\label{K0Lob}
K_0(t)=\left\{\begin{array}{ll}
-1-t+\lambda,&-1\le t\le -x,\qquad\qquad\qquad\quad\\[2mm]
-t,&-x<t\le x,\\[2mm]
\ \ 1-t-\lambda,&\ \ \,x<t\le 1,	
\end{array}\right.
\end{equation}
and
\begin{equation}\label{K1Lob}
K_1(t)=\left\{\begin{array}{ll}
\dfrac{1}{2}(1+t)(1+t-2\lambda),&-1\le t\le -x,\\[3mm]
\dfrac{1}{2}(1 + t^2 - 2 x(1-\lambda)-2\lambda),&-x<t\le x,\\[3mm]
\dfrac{1}{2}(1 - t)(1-t-2\lambda),&\ \ \,x<t\le 1,	
\end{array}\right.
\end{equation}
and presented in Fig.~\ref{AlomFig1} (left).
Since 
\begin{equation}\label{M0eq34}
M_0(x)=\int_{-1}^1|K_0(t)|{\D}t=\left\{\begin{array}{ll}
\lambda^2+x^2+(1-\lambda-x)^2,&0\le x\le 1-\lambda,\\[2mm]
2x(1-\lambda)+2\lambda-1,&1-\lambda<x\le 1,	
\end{array}\right.
\end{equation}
the corresponding inequality of Ostrowski type is
\begin{equation}\label{m0op}
\left|\int_{-1}^1f(t){\D}t- Q_{4,0}^\lambda(f)\right|\le 
M_0(x)\|f'\|_\infty,\quad 0\le x\le 1.
\end{equation}
This result for $x\in[0,1-\lambda]$ has been obtained by Alomari \cite{Alom12}, including several particular cases. 

The optimal estimate (\ref{m0op}) is attained when
\[\min_{x\in[0,1]}M_0(x)=M_0\Bigl(\frac{1-\lambda}{2}\Bigr)=\frac{1}{2} \left(3 \lambda ^2-2 \lambda +1\right).\]
Moreover, this value becomes the smallest ($1/3$) for $\lambda=1/3$, when the rule (\ref{Alom12}) reduces to
\[ Q_{4,0}^{1/3}(f)= \frac{1}{3}[f(-1)+f(1)]+\frac{2}{3}\left[ f\Bigl(-\frac{1}{3}\Bigr)+f\Bigl(\frac{1}{3}\Bigr)\right].
\]

As we mentioned before, for $\lambda=0$ and $\lambda=1$ the rule (\ref{Alom12}) reduces to the case considered by Guessab and Schmeisser \cite{GSch}. For $x=0$, 
the rule (\ref{Alom12}) becomes a symmetric three-point rule 
\[ Q_{3,0}^\lambda(f)= \lambda[f(-1)+f(1)]+2(1-\lambda)f(0),\]
which is a special case of (\ref{DCR2000}) (for $x=0$).

The corresponding result for the bound $M_1(x)$, for $r=1$ and $1/2\le \lambda<1$,  can be found in the form $M_1(x)=(1-\lambda)x^2+\lambda-\frac{1}{3}$. However, for $0<\lambda<1/2$, the expression for this bound is quite complicated. For example, for 
$\lambda=1/5$ we have
\[M_1(x)=\int_{-1}^1|K_1(t)|{\D}t=\left\{ 
\begin{array}{ll}
\dfrac{58}{375}-\dfrac{4 x^2}{5},&0<x\le \dfrac{3}{8},\\[3mm]
\dfrac{2}{375}\left(29-150 x^2+10 \sqrt{5} (8x-3)^{3/2}\right),&
\dfrac{3}{8}<x<\dfrac{3}{5},\\[3mm]
\dfrac{2}{15} \left(6 x^2-1\right),&\dfrac{3}{5}\le x<1.	
\end{array}\right.\] 
The bound $M_1(x)$, $0\le x\le 1$, in the inequality of Ostrowski's type
\[\left|\int_{-1}^1f(t){\D}t-Q_{4,0}^\lambda(f)\right|\le 
M_1(x)\|f''\|_\infty,\]
for different values of the parameter $\lambda$ are displayed in Fig.~\ref{AlomFig1} (right).
 
In the sequel of this subsection we consider the case of (\ref{Alom12}), when its degree of precision is $d=5$, i.e., when
$\lambda=1/6$ and  $x=\xi=1/\sqrt{5}$. In fact, it is a Lobatto quadrature formula with two internal nodes (cf. \cite[pp.~330--332]{GMGVM08}). In order to get  the following estimates
\[\left|\int_{-1}^1f(t){\D}t- \frac{1}{6}\bigl[f(-1)+f(1)\bigr]-\frac{5}{6} \left[f\left(-\frac{1}{\sqrt{5}}\right)+f\left(\frac{1}{\sqrt{5}}\right)\right]\right|\le M_r\|f^{(r+1)}\|_\infty,\]
for each $r=0,1,\ldots,d\,(=5)$, where $\DS M_r=\int_{-1}^1|K_r(t)|{\D}t$, we need the corresponding kernels.\\ 

For $r=0$ and $r=1$, the  expressions for $K_r(t)$ are given by (\ref{K0Lob}) and (\ref{K1Lob}) for arbitrary $\lambda$ and $x$, and we must take  $\lambda=1/6$ and  $x=1/\sqrt{5}$, while  for $2\le r\le 5$ we have the following expressions
\begin{figure}[ht]
\centering
\includegraphics[width=0.48\textwidth]{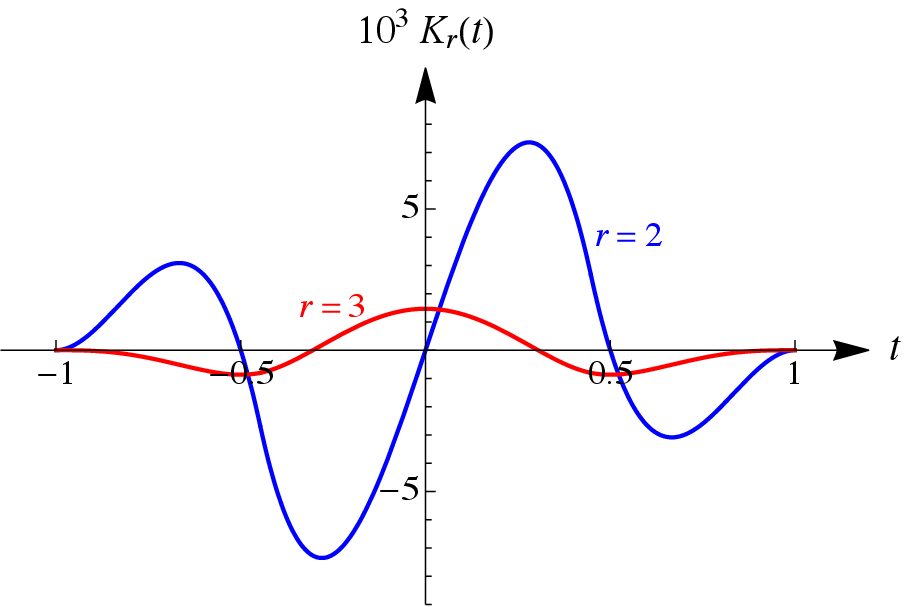}\quad
\includegraphics[width=0.48\textwidth]{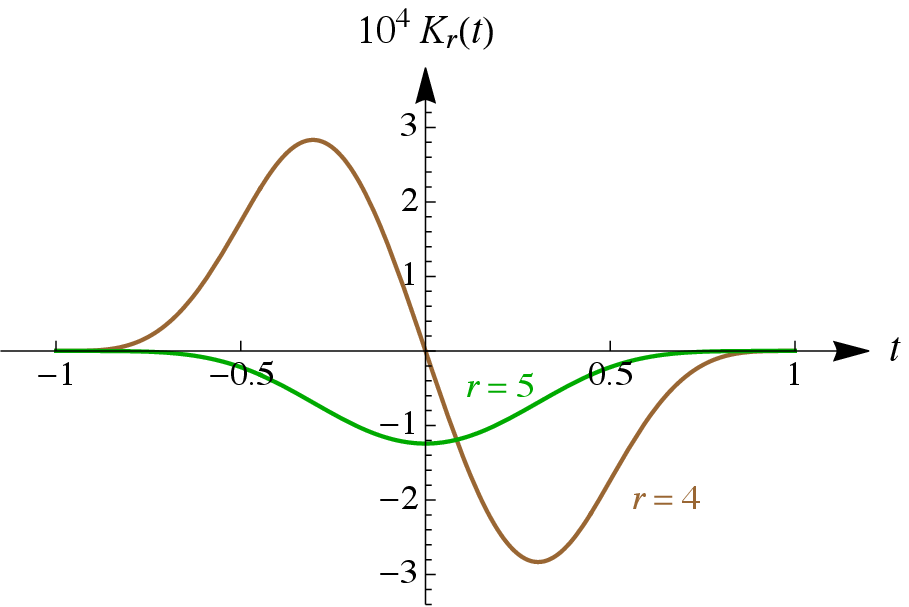}
\caption[]{The  kernels $K_2(t)$ and $K_3(t)$  (left)  and 
$K_4(t)$ and $K_5(t)$ (right) for the Lobatto quadrature rule  with two internal nodes}\label{AlomFig2}	
\end{figure}
\[
K_2(t)=\left\{\begin{array}{ll}
-\dfrac{1}{12}(1+t)^2(1+2t),&\ \, -1\le t\le -\dfrac{1}{\sqrt{5}},\ \,\qquad\qquad\qquad\qquad\qquad\\[3mm]
\,-\dfrac{1}{6}t(2-\sqrt{5}+ t^2),&-\dfrac{1}{\sqrt{5}}<t\le \dfrac{1}{\sqrt{5}},\ \\[3mm]
\ \ \dfrac{1}{12}(1 - t)^2(1-2t),&\ \ \,\dfrac{1}{\sqrt{5}}<t\le 1,	
\end{array}\right.
\]
\[
K_3(t)=\left\{\begin{array}{ll}
 \dfrac{1}{72}(1+t)^3(1+3t),&\ \, -1\le t\le -\dfrac{1}{\sqrt{5}},\qquad\quad\\[3mm]
 \dfrac{1}{360}\left(5-2\sqrt{5}+ 30(2-\sqrt{5})t^2+15t^4\right),&-\dfrac{1}{\sqrt{5}}<t\le \dfrac{1}{\sqrt{5}},\\[3mm]
 \dfrac{1}{72}(1 - t)^3(1-3t),&\ \ \,\dfrac{1}{\sqrt{5}}<t\le 1,	
\end{array}\right.
\]
\[
K_4(t)=\left\{\begin{array}{ll}
 -\dfrac{1}{720}(1+t)^4(1+6t),&\ \, -1\le t\le -\dfrac{1}{\sqrt{5}},\qquad\quad\\[3mm]
\ \ \dfrac{1}{360}t\left(2\sqrt{5}-5+ 10(\sqrt{5}-2)t^2-3t^4\right),&-\dfrac{1}{\sqrt{5}}<t\le \dfrac{1}{\sqrt{5}},\\[3mm]
\ \ \dfrac{1}{720}(1 - t)^4(1-6t),&\ \dfrac{1}{\sqrt{5}}<t\le	1,\end{array}\right.
\]
\[
K_5(t)=\left\{\begin{array}{ll}
 \ \dfrac{1}{720}t(1+t)^5,&\ \, -1\le t\le -\dfrac{1}{\sqrt{5}},\\[3mm]
\ \dfrac{1}{720}\biggl[t^6-5(\sqrt{5}-2)t^4+(5-2\sqrt{5})t^2-\dfrac{\sqrt{5}}{25}\,\biggr],&-\dfrac{1}{\sqrt{5}}<t\le \dfrac{1}{\sqrt{5}},\\[3mm]
\!\!-\dfrac{1}{720}t(1 - t)^5,&\ \ \,\dfrac{1}{\sqrt{5}}<t\le	1.\end{array}\right.
\]
Their graphics are presented in Fig.~\ref{AlomFig2}. 

Finally, we get values for the bounds $M_r$, $0\le r\le 5$,
\begin{align*} 
M_0=&\frac{1}{90}\left(101-30 \sqrt{5}\right)\approx 0.376866,\\[2mm]  
M_1=&\frac{1}{81} \left(1+12 \sqrt{3 \left(17 \sqrt{5}-38\right)}\right)\approx 0.04177718,\qquad\qquad\qquad\\[2mm] 
M_2=&\frac{45-16 \sqrt{5}}{1440}\approx 6.40480025\times 10^{-3},\\[2mm] 
M_3=&1.132646548 \times 10^{-3},\\[2mm]
M_4=&\frac{1}{1800 \sqrt{5}}\approx 2.48451997\times 10^{-4},\\[2mm] 
M_5=&\frac{2}{23625}\approx 8.4656084656\times 10^{-5}.\qquad\qquad\qquad\qquad\qquad\qquad
\end{align*}
The same values of these constants can be found in the book 
\cite[\S4.2.1]{Mon2012}.   

\vspace{0.8cc}

\noindent {\large\bf 3.8. Four point  inequality with double internal nodes by Liu and Park  \cite{LP}}
\vspace{0.5cc} 

There are many inequalities of Ostrowski's type, with including derivatives in quadrature sums. Here, we consider a symmetric four-point quadrature rule with double internal nodes \cite{LP}  
\begin{equation}\label{LiuPark}
Q_{4,2}(f)=\frac{1}{2}\left[f(-1)+f(-x)+f(x)+f(1)\right]-\frac{x}{2}\left[f'(x)-f'(-x)\right].	
\end{equation}

\noindent
{\bf Remark.} A three-point formula with a double inner node
i.e.,
\[Q_{3,1}(f)=\frac{1}{2}\bigl[f(-1)+2f(x)+f(1)\bigr]-x f'(x),\]
where $-1\le x\le1$, was considered by Dragomir and Sofo \cite{DS2000}. 
\vspace{0.5cc} 

Without loss of generality we suppose that $0\le x\le1$ in the rule (\ref{LiuPark}). This rule  has degree of exactness $d=1$, except $x=1/\sqrt{3}$, because of
\begin{equation}\label{stex}
\left\{R_{4,2}(e_r)\right\}_{r=0}^\infty=
\left\{0,\,0,\,x^2-\frac{1}{3},\,0,\,3 x^4-\frac{3}{5},\,0,\,\ldots\right\}.	
\end{equation}
Since
\[K_0(t)=\left\{\begin{array}{ll}
-\dfrac{1}{2}(1+2t),&-1<t\le -x,\\[3mm]
-t,&-x<t\le x,\\[3mm]
\ \ \dfrac{1}{2}(1+2t),&\ \ x<t\le 1, 	
\end{array}\right.\]
and
\[K_1(t)=\left\{\begin{array}{ll}
\dfrac{1}{2}t(1+t),&-1<t\le -x,\\[4mm]
\dfrac{1}{2}t^2,&-x<t\le x,\\[4mm]
\dfrac{1}{2}t(-1+t),&\ \ x<t\le 1, 	
\end{array}\right.\]
with
\[M_0(x)=\int_{-1}^1|K_0(t)|{\D}t=\left\{\begin{array}{ll}
\dfrac{1}{2}-x+2x^2,&0<x<\dfrac{1}{2},\\[4mm]
x,&x\ge \dfrac{1}{2}. 	
\end{array}\right.\]
and
\[M_1(x)=\int_{-1}^1|K_1(t)|{\D}t=\frac{1}{6}(1-3x^2+4x^3),\]
respectively, the corresponding inequalities of Ostrowski's type are
\[\left|I(f)-Q_{4,2}(f)\right|\le M_r(x)\|f^{(r+1)}\|_\infty\quad (r=0,1).\]
The minimal values of $x\mapsto M_0(x)$ and $x\mapsto M_1(x)$ are
\[\min_{x\in[0,1]}M_0(x)=M_0\Bigl(\frac{1}{4}\Bigr)=\frac{3}{8}\quad\mbox{and}\quad
\min_{x\in[0,1]}M_1(x)=M_1\Bigl(\frac{1}{2}\Bigr)=\frac{1}{8}.\]
The case $r=1$ is considered in \cite[Theorem 2.1]{LP}.

For $x=1/\sqrt{3}$ the rule (\ref{LiuPark}) reduces to
\begin{align}\label{haQ42}
\widehat Q_{4,2}(f)=&\frac{1}{2}\left[f(-1)+f\Bigl(-\dfrac{1}{\sqrt{3}}\Bigr)+f\left(\dfrac{1}{\sqrt{3}}\right)+f(1)\right]\\[2mm]
&\qquad\qquad\quad -\dfrac{1}{2\sqrt{3}}\left[f'\Bigl(\dfrac{1}{\sqrt{3}}\Bigr)-f'\Bigl(-\dfrac{1}{\sqrt{3}}\Bigr)\right],\nonumber
\end{align}
and according to (\ref{stex}), its  degree of exactness is $d=3$. 
\begin{figure}[h]
\centering
\includegraphics[width=0.65\textwidth]{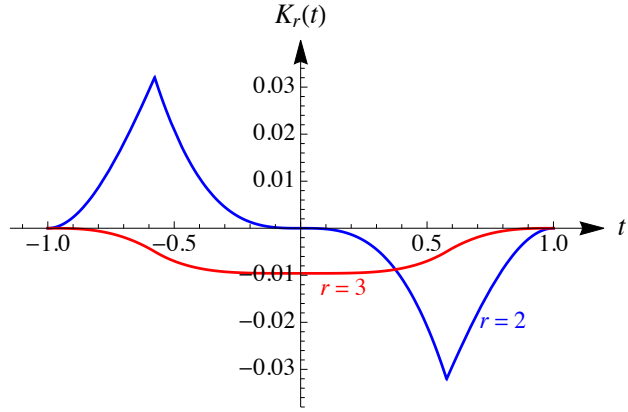} 
\caption[]{The kernels  $K_2(t)$ and $K_3(t)$}
\label{sl1}	
\end{figure}

For $r=2$ and $r=3$ the corresponding Peano kernels are (see Fig.~\ref{sl1})
\[K_2(t)=\left\{\begin{array}{ll}
\ \, \dfrac{1}{12}(1+t)^2(1-2t),&\ -1<t\le -\dfrac{1}{\sqrt{3}},\\[4mm]
-\dfrac{1}{6}t^3,&-\dfrac{1}{\sqrt{3}}<t\le \dfrac{1}{\sqrt{3}},\\[4mm]
-\dfrac{1}{12}(1-t)^2(1+2t),&\ \ \dfrac{1}{\sqrt{3}}<t\le 1, 	
\end{array}\right.\]
and
\[K_3(t)=\left\{\begin{array}{ll}
\dfrac{1}{24} (t-1) (t+1)^3,&-1<t\le -\dfrac{1}{\sqrt{3}},\\[4mm]
\dfrac{1}{216} \left(9 t^4+4 \sqrt{3}-9\right),&-\dfrac{1}{\sqrt{3}}<t\le \dfrac{1}{\sqrt{3}},\\[4mm]
\dfrac{1}{24} (t-1)^3 (t+1),&\ \ \dfrac{1}{\sqrt{3}}<t\le 1. 	
\end{array}\right.\]

Since
\[M_2=\int_{-1}^1|K_2(t)|{\D}t=\frac{9-4\sqrt{3}}{108}\quad\mbox{and}\quad  M_3=\int_{-1}^1|K_3(t)|{\D}t=\frac{1}{90},\]
we get the following inequalities of Ostrowski type are
\[\left|I(f)-\widehat Q_{4,2}(f)\right|\le M_r \|f^{(r+1)}\|_\infty \quad (r=2,3),\]
where $\widehat Q_{4,2}(f)$ is given by (\ref{haQ42}).

\vspace{1.cc}
\begin{center}
{\bf 4. CONCLUSION}
\end{center}

In this survey paper we considered only selected simple inequalities of Ostrowski's type and their estimates for differentiable functions with bounded derivative. In fact, our examples are special cases of an general four-point quadrature formula with double nodes
\begin{align*} 
Q_{4,4}(f)=&\lambda[f(1)+f(-1)]+(1-\lambda)[f(x)+f(-x)]\\[2mm]
&\qquad\ \ +\gamma[f'(1)-f'(-1)]+\delta[f'(x)-f'(-x)],\nonumber
\end{align*}
where $\lambda,\gamma,\delta$ are real parameters,  $0<\lambda<1$ and $0< x<1$. It could be interesting to investigate this general case and analyse its  particular cases. Finally, we mention some papers which deal with the weighted inequalities \cite{AgPeTi,KoPeTi2013,MaPeUj2000,PeRiVu2008,AW,UjLe2008}, as well as the book \cite{Mon2012}.

\vspace{1cc}


{\small
\noindent
Serbian Academy of Sciences and Arts, Beograd, Serbia\\
\&\\
University of Ni\v s, Faculty of Sciences and Mathematics, 
Ni\v s, Serbia\\
E-mail: {\tt gvm@mi.sanu.ac.rs}} 

\vspace{1cc}


\end{document}